\magnification\magstep1

\def\n{\noindent}
\def\vp{\varepsilon}
\def\Bbb#1{{\bf#1}}
\def\bb#1{{\bf#1}}
\def\ca#1{{\cal#1}}
\def\brel#1{{\buildrel #1 \over \longrightarrow}}
\def\brell#1{{\buildrel #1 \over {\hbox to 35pt{\rightarrowfill}}}}

\overfullrule = 0pt
\settabs 5\columns
\pageno=0

\centerline{FACTORIZATION OF COMPLETELY BOUNDED}
\centerline{BILINEAR OPERATORS AND INJECTIVITY}\vskip.5in

\+Allan M.\ Sinclair$^{(*)}$&&&Roger R.\ Smith$^{(*)(\dag)}$\cr
\+Department of Mathematics&&&Department of Mathematics\cr
\+University of Edinburgh&&&Texas A\&M University\cr
\+Edinburgh EH9 3JZ&&&College Station, TX \ 77843\cr
\+SCOTLAND&&&U.S.A.\cr
\+allan@mathematics.edinburgh.ac.uk&&&rsmith@math.tamu.edu\cr
\vskip1.5in

\centerline{ABSTRACT}\medskip

A completely bounded bilinear operator $\phi\colon \ \ca M\times \ca M\to
\ca M$ on a von~Neumann algebra $\ca M$ is said to have a factorization in
${\cal M}$ if
there  exist completely bounded linear operators\break $\psi_j$,
$\theta_j\colon \ \ca M\to \ca M$ such that
$$\phi(x,y) = \sum_{j\in\Lambda} \psi_j(x) \theta_j(y),\qquad x,y\in \ca
M,$$
where convergence of the sum is made precise below. The main result of the
paper is that all completely bounded bilinear operators $\phi\colon \ \ca M
\times \ca M\to \ca M$ have factorizations in $\ca M$ if and only if $\ca
M$ is injective.

\vfil

\n $^{(*)}$Partially supported by a NATO collaborative research grant.

\n $^{(\dag)}$Partially supported by an NSF research grant.\eject

\baselineskip = 18pt

{\bf \S 1. Introduction}

There are several conditions on a von~Neumann algebra $\ca N$ that are
known to be equivalent to the injectivity of $\ca N$. The outstanding, and
fundamental, result is Connes' proof [10] that injective factors on a
separable Hilbert space are hyperfinite (see also [32]). Subsequently
Haagerup [19] and Popa [26] gave simpler treatments of this result which
avoided the technical theory of automorphism groups of von~Neumann algebras
in [10]. One result in the development of the subject prior to [10] plays a
role below. Effros and Lance [15, Corollary 4.6] showed that a von~Neumann
factor $\ca N$ is semidiscrete (equivalently injective) if and only if the
$C^*$-algebra $C^*(\ca N, \ca N')$ is isomorphic to $\ca N\otimes_{\rm min}
\ca N'$; this is used in proving Theorem~4.4 below. The operators between
von~Neumann
algebras which appear in [15] are all completely positive, but there are
characterizations of injectivity of a von~Neumann algebra $\ca N$ based on
properties of completely bounded linear operators associated with $\ca N$.
For example, Haagerup [20] has shown that $\ca N$ is injective if and only
if each completely bounded linear operator from $\ell^\infty$ into $\ca N$
is a linear combination of completely positive linear operators from
$\ell^\infty$ into $\ca N$; this is used in proving  Theorem~5.3 below.

Let $\ca M$ and $\ca N$ be von~Neumann algebras with $\ca N$  acting on a
Hilbert space $H$ and $\ca M$ infinite dimensional. The representation
theorem of a completely bounded bilinear operator from $\ca M\times\ca M$
into $B(H)$ provides a factorization of such an operator into $\ca N$.
Strengthening the hypotheses on this factorization for all completely
bounded bilinear operators $\phi\colon \ \ca M\times \ca M\to\ca N$ provides
another characterization of injectivity of $\ca N$ as we shall explain. If
$\phi\colon \ \ca M\times \ca M \to \ca N\subseteq B(H)$ is a completely
bounded bilinear operator then there is a representation $\pi\colon \ \ca M
\to B(K)$ and continuous linear operators $W\colon \ H\to K$, $T\colon \
K\to K$, and $V\colon \ K\to H$ such that
$$\phi(m_1,m_2) = V\pi(m_1) T\pi(m_2)W,\qquad m_1,m_2\in \ca M\eqno (1.1)$$
and $\|\phi\|_{cb} = \|V\|\, \|T\|\, \|W\|$. (See [5, 6, 23, 27]).
However there is little control over $V,T$, and $W$ other than the norm
estimate. If $\psi,\theta\colon \ \ca M\to \ca N$ are completely bounded
then $\phi\colon \ \ca M\times \ca M\to \ca N $ defined by
$$\phi(m_1,m_2) = \psi(m_1) \theta(m_2),\qquad m_1,m_2\in \ca M\eqno
(1.2)$$
is a completely bounded bilinear operator, and (1.2) represents a
factorization of the bilinear operator $\phi$ as a product of linear
operators. More generally, suitable weakly convergent sums
$\sum\limits_{j\in\Lambda} \psi_j(m_1) \theta_j(m_2)$ of such products may
define a completely bounded bilinear operator $\phi\colon \ \ca M\times \ca
M\to \ca N$, and we refer to such a sum as a factorization of $\phi$. The
main result of the paper is that injectivity of $\ca N$ is equivalent to
all completely bounded operators $\phi\colon \ \ca M\times \ca M\to \ca N$
having such factorizations. A consequence of our work is that
factorizations with $\psi_j,\theta_j$ mapping into $B(H)$ are always
possible; the crucial point is to require $\psi_j,\theta_j$ to map into
$\ca N$. Indeed it suffices to take $\ca M = \ca N$, which gives a
characterization of injectivity in terms of completely bounded bilinear
operators which is internal to $\ca N$.

We now give a brief description of the contents of the paper. Section~2
contains the basic notation and definitions, and also a short account of
the $w^*$-Haagerup tensor product of $CB(\ca X,\ca N) \otimes_{w^*h} CB(\ca
Y,\ca N)$, where $\ca X, \ca Y$ are operator spaces, $\ca N$ is a
von~Neumann algebra, and $CB(\ca X,\ca N)$ is the space of completely
bounded linear operators from $\ca X$ into $\ca N$. This tensor product
provides a convenient language for the formulation of our results. However
we have delayed its appearance until the last section to reduce the
technicalities for readers who are unfamiliar with it.

Section 3 contains a theorem on module map extensions in the bilinear case,
extending a result of Wittstock [33] for one variable. This is used to
obtain Proposition~3.3, a technical result on the representation of modular
bilinear operators, which is important subsequently. The fourth and fifth
sections are the heart of the paper, each devoted to one implication in the
equivalence of injectivity and the factorization of completely bounded
bilinear operators. Injectivity implies factorization is Theorem~4.4, while
the reverse implication is Theorem~5.4. The final section is a brief
summary of results and includes some other equivalences formulated in terms
of the $w^*$-Haagerup tensor product.

We refer the reader to [22] for an account of the theory of completed
bounded linear operators. The subsequent development of the multilinear
case may be found in the survey article [6] or the book [27]. We also
refer to [3, 7, 8, 19, 21, 24, 25, 26, 32] for related results on
injectivity and multilinear operators.
\vfill\eject

\n {\bf \S 2. Notations and definitions}

Throughout $\ca M$ and $\ca N$ will denote von~Neumann algebras acting on a
Hilbert space $H$ with commutants $\ca M'$ and $\ca N'$. $C^*$-algebras are
denoted by $\ca A$ and operator spaces by $\ca E, \ca F$, $\ca X$ or $\ca
Y$. Recall that an operator space $\ca X$ is a norm closed subspace of
$B(H)$, the algebra of bounded linear operators on $H$, together with the
norms and structure of $\bb M_n(\ca X)$ in $\bb M_n(B(H)) = B(H^n)$, where
$\bb M_n$ denotes the $n\times n$ matrices. We refer to [22] and [6, 27]
respectively for the theories of completely bounded linear operators and
completely bounded multilinear operators. Recall that a completely bounded
linear operator $\phi\colon \ \ca X\to B(H)$, where $\ca X$ is an operator
space in a $C^*$-algebra $\ca A$, has a representation of the form
$$\phi(x) = U\pi(x)V,\qquad x\in \ca X.\eqno (2.1)$$
Here $\pi$ is a representation of $\ca A$ on a Hilbert space $K$ and
$V\colon
\ H\to K$, $U\colon \ K\to H$ are continuous linear operators satisfying
$$\|U\| = \|V\| = \|\phi\|^{1/2}_{cb}.\eqno (2.2)$$
The corresponding result for completely bounded bilinear operators is the
following. Given operator spaces $\ca X$ and $\ca Y$ in a $C^*$-algebra
$\ca A$ and a completely bounded bilinear operator $\phi\colon \ \ca
X\times \ca Y\to B(H)$, there exist a representation $\pi\colon \ \ca A\to
B(K)$ and continuous linear operators $W\colon \ H\to K$, $T\colon \ K\to
K$, and $U\colon \ K\to H$ such that
$$\phi(x,y) = U\pi(x) T\pi(y)W,\qquad x\in \ca X,y \in \ca Y.\eqno (2.3)$$
Moreover $U,T$ and $W$ may be  chosen to satisfy the optimal condition
$$\|U\| = \|T\| = \|W\| = \|\phi\|^{1/3}_{cb}.\eqno (2.4)$$
The usual formulation of (2.3) is
$$\phi(x,y) = U\pi(x) T\rho(y)W,\qquad x\in \ca X, y\in \ca Y\eqno (2.5)$$
where $\pi$ and $\rho$ are possibly distinct representations of $\ca A$,
but this may be reduced to the form of (2.3) by writing
$$\phi(x,y) = (U,0) \left(\matrix{\pi(x)&0\cr 0&\rho(x)\cr}\right)
\left(\matrix{0&T\cr 0&0\cr}\right) \left(\matrix{\pi(y)&0\cr 0&\rho(y)\cr}
\right) {0\choose W}.\eqno (2.6)$$

If $\ca X, \ca Y$ and $\ca Z$ are operator spaces, $CB(\ca X,\ca Y)$
denotes the space of completely bounded linear operators of $\ca X$ into
$\ca Y$, while $CB^2(\ca X\times \ca Y,\ca Z)$ denotes the space of
completely bounded bilinear operators of $\ca X\times \ca Y$ into $\ca Z$.
When $\ca X= \ca Y$, we abbreviate this latter space to $CB^2(\ca X,\ca
Z)$.

For $\phi \in CB^2(\ca A, B(H))$, recall that the $n$-fold amplification
$\phi_n\in CB^2(\bb M_n(\ca A), \bb M_n(B(H)))$ is defined by
$$\phi_n((x_{ij}), (y_{ij})) = \left(\sum^n_{k=1} \phi(x_{ik},
y_{kj})\right) \eqno (2.7)$$
for $(x_{ij}), (y_{ij}) \in \bb M_n(\ca A)$. Then $\phi$ is said to be
completely positive if
$$\phi_n((x_{ij}), (x_{ij})^*)\ge 0,\quad (x_{ij}) \in \bb M_n(\ca A),
\qquad n\in\bb N.\eqno (2.8)$$
In contrast to the linear case, completely positive bilinear operators need
not be completely bounded. This is well known [5], but we include an
elementary example. Let $\psi\colon \ B(\ell_2)\to B(\ell_2)$ be the
transpose on infinite matrices	and define $\phi\colon \ B(H) \times
B(H)\to B(H)$ by
$$\phi(x,y) = \psi(x) \psi(y^*)^*,\qquad x,y\in B(H).\eqno (2.8)$$
It is easy to check that
$$\phi_n((x_{ij}), (x_{ij})^*) = \psi_n(x_{ij}) \psi_n(x_{ij})^* \ge
0,\eqno (2.9)$$
and so $\phi$ is completely positive. However $\psi(x) = \phi(x,1)$, and so
$\phi$ cannot be completely bounded, since $\psi$ is not. Thus decomposing
bilinear completely bounded operators as a linear combination of completely
positive bilinear operators is (seemingly) less restrictive than similar
decompositions in the linear case. We note in passing that the representation
(2.3) may be extended to the multilinear case [5, 23], but will not be
needed here.

Let $\ca X, \ca Y$ and $\ca Z$ be bimodules over a $C^*$-algebra $\ca A$
and
let $\phi\colon \ \ca X\times \ca Y\to \ca Z$ be bilinear. Then $\phi$ is
said to be $\ca A$-modular if, for $x\in \ca X$, $y\in \ca Y$, $a\in\ca A$,
the following relations hold:
$$\eqalignno{\phi(ax,y) &= a\phi(x,y),&(2.10)\cr
\phi(xa,y) &= \phi(x,ay),&(2.11)\cr
\phi(x,ya) &= \phi(x,y)a.&(2.12)}$$
Such operators have played a major role in the theory of completely bounded
operators and their applications for several years (see [4, 8, 13, 27,
28]), and will also be useful in subsequent sections of this paper. Recall
that a von~Neumann algebra $\ca N$ is injective if for any containing
von~Neumann algebra $\ca M$ there is a conditional expectation $\bb E\colon
\ \ca M\to \ca N$. By this we mean a completely positive projection of $\ca
M$ onto $\ca N$, and such projections are automatically $\ca N$-modular
[31].

The $C^*$-algebra generated by two $C^*$-subalgebras $\ca A$ and $\ca B$ of
$B(H)$ is denoted by $C^*(\ca A,\ca B)$. If $\Lambda$ is a (non-empty)
index set and $H$ is a Hilbert space, let $\ell_2(\Lambda, H)
=\ell_2(\Lambda)\otimes_2 H$ denote the Hilbert space of ``sequences'' in
$H$ indexed by $\Lambda$. For a minimal projection $e$ onto a standard basis
vector of $\ell_2(\Lambda)$, we let
$$\eqalign{R(\Lambda) &= B(\ell_2(\Lambda))e,\cr
C(\Lambda) &= eB(\ell_2(\Lambda)),}$$
be the $\Lambda$-row and -column operator spaces respectively. Although we
require general index sets for the proper formulation of our results, the
reader will not be misled by thinking of $\Lambda$ as $\bb N$.

We now review several tensor products which will be needed subsequently.
The minimal (also called injective or spatial) tensor product of
$C^*$-algebras $\ca A$ and $\ca B$ is denoted by $\ca A\otimes_{\rm min}
\ca
B$ [30], while $\ca M\overline\otimes \ca N$ denotes the von~Neumann
algebra tensor product of von~Neumann algebras $\ca M$ and $\ca N$ [30].
The Haagerup
tensor product $\ca A\otimes_h \ca B$ [13, 18] is the completion  of the
algebraic tensor product $\ca A\otimes\ca B$ in the norm
$$\|u\|_h = \inf\left\{\left\|\sum^n_{j=1} a_ja^*_j\right\|^{1/2}
\left\|\sum^n_{j=1} b^*_jb_j\right\|^{1/2}\right\}\eqno (2.13)$$
taken over all representations $u = \sum\limits^n_{j=1} a_j \otimes b_j \in
\ca A \otimes \ca B$. There are several weak versions of this tensor
product and we will require the $w^*$-Haagerup tensor product
$\otimes_{w^*h}$, introduced for pairs of dual operator spaces in [2]. Our
interest will focus on $CB(\ca X, \ca N) \otimes_{w^*h} CB(\ca Y,\ca N)$
where $\ca X$ and $\ca Y$ are operator spaces and $\ca N\subseteq B(H)$ is
a von~Neumann algebra, and we give a straightforward definition in this
case which is equivalent to the original formulation. Note that $CB(\ca
X,\ca N)$ is a dual operator space which can be identified with the dual of
the operator space projective tensor product $\ca X\widehat\otimes \ca N_*$
[1, 16]. We omit further discussion of $\widehat\otimes$ since it will
not be needed subsequently.

Consider the vector space $\ca V$ of all formal sums
$\sum\limits_{j\in\Lambda} \psi_j \otimes \theta_j$ where $\psi_j \in
CB(\ca X,\ca N)$, $\theta_j \in CB(\ca Y,\ca N)$ and, for all $x\in \ca X$,
$y\in \ca Y$, and finite subsets $F$ of $\Lambda$, there exists a constant
$K$ such that
$$\eqalignno{\left\|\sum_{j\in F} \psi_j(x) \psi_j(x)^*\right\| &\le
K\|x\|^2, &(2.14)\cr
\left\|\sum_{j\in F} \theta_j(y)^* \theta_j(y)\right\| &\le
K\|y\|^2.&(2.15)}$$
For vectors $\xi,\eta \in H$,
$$\eqalignno{\sum_{j\in F} \|\psi_j(x)^*\eta\|^2 &= \sum_{j\in F} \langle
\psi_j(x) \psi_j(x)^* \eta,\eta\rangle\cr
&\le K\|x\|^2 \|\eta\|^2,&(2.16)}$$
and similarly
$$\sum_{j\in F} \|\theta_j(y)\xi\|^2 \le K\|y\|^2 \|\xi\|^2.\eqno (2.17)$$
These two inequalities then remain valid for the sum over all $j\in
\Lambda$ (as a consequence of which only countably many terms in the sums
are non-zero), and the Cauchy-Schwarz inequality then shows that
$\sum\limits_{j\in\Lambda} \langle\psi_j(x) \theta_j(y)\xi,\eta\rangle$ is
an absolutely convergent series, bounded in absolute value by $K\|x\|\,
\|y\|\, \|\xi\|\, \|\eta\|$. Thus $\sum\limits_{j\in\Lambda} \psi_j(x)
\theta_j(y)$ is weakly convergent, so defines an element of $\ca N$.
Essentially the same argument shows that
$$\sum_{j\in\Lambda} \langle \psi_j(x) t\theta_j(y) \xi,\eta\rangle, \qquad
x\in \ca X, y\in \ca Y, t\in B(H), \xi,\eta\in H\eqno (2.18)$$
is always an absolutely convergent series and so $\sum\limits_{j\in\Lambda}
\psi_j(x) t \theta_j(y)$ converges weakly to an element of $B(H)$ for all
$t\in B(H)$. We now declare two sums $\sum\limits_{j\in\Lambda}
\psi_j\otimes \theta_j$, $\sum \widetilde\psi_j \otimes \tilde\theta_j$ to
be equal in $\ca V$ if
$$\sum_{j\in \Lambda} \psi_j(x) t\theta_j(y) = \sum_{j\in\Lambda}
\widetilde \psi_j(x)t\tilde \theta_j(y)\eqno (2.19)$$
for all $t\in B(H)$. By taking $t=I$ it is then clear that
$\sum\limits_{j\in\Lambda} \psi_j(x) \theta_j(y)$ is independent of the
particular representation chosen.

An element $v = \sum\limits_{j\in\Lambda} \psi_j\otimes \theta_j$ leads to
bounded bilinear maps $\Psi\colon \ \ca X\times \ca X^*\to \ca N$ and
$\Theta\colon \ \ca Y^*\times \ca Y \to \ca N$ defined by
$$\eqalignno{\Psi(x_1,x^*_2) &= \sum_{j\in\Lambda} \psi_j(x_1)
\psi_j(x_2)^*,&(2.20)\cr
\Theta(y^*_1,y_2) &= \sum_{j\in\Lambda} \theta_j(y_1)^*
\theta_j(y_2),&(2.21)}$$
with weak convergence in the sums.  For each $v = \sum\limits_{j\in\Lambda}
\psi_j \otimes \theta_j\in {\cal V}$ we define
$$|||v||| = \inf\{\|\Psi\|^{1/2}_{cb} \|\Theta\|^{1/2}_{cb}\}\eqno (2.22)$$
where the infimum is taken over all possible representations of $v$. If no
representation of $v$ has associated operators which are completely bounded,
then we set $|||v||| = \infty$. It is not immediately clear that
$|||\cdot|||$ is a norm on the set of elements for which $|||v|||$ is
 finite. This will follow
from the next proposition, whose purpose is to give an alternative
description of the $w^*$-Haagerup norm on $CB({\cal X}, {\cal N})
\otimes_{w^*h} CB({\cal Y},{\cal N})$.

\proclaim Proposition 2.1. Let $v\in {\cal V}$. Then $v\in CB({\cal X},
{\cal N})\otimes_{w^*h} CB({\cal Y}, {\cal N})$ if and only if $|||v||| <
\infty$, and in this case
$$|||v||| = \|v\|_{w^*h}.\eqno (2.23)$$

\n {\bf Proof.} Suppose that $|||v|||  < \infty$. Since it is clear that
$|||\lambda v||| = |\lambda|\ |||v|||$, we may assume that
$|||v|||=1$. Thus, given $\vp>0$, there exists a representation $v =
\sum\limits_{j\in\Lambda} \psi_j \otimes\theta_j$ and $\|\Psi\|_{cb}$,
$\|\Theta\|_{cb} < (1+\vp)^{1/2}$. For any finite subset ${\cal F} =
\{1,\ldots, n\} \subseteq \Lambda$ (after renumbering), the norm of
$(\psi_1,\ldots, \psi_n)$ as an element of ${\Bbb M}_n(CB({\cal X}, {\cal
N}))$ is its $cb$-norm as an element of $CB({\cal X}, {\Bbb M}_n({\cal
N}))$, so, for $X\in {\Bbb M}_k({\cal X})$, $k\ge 1$,
$$\eqalignno{&~\|(\psi_1\otimes I_k(X),\ldots, \psi_n\otimes I_k(X))\|^2\cr
= &\left\|\sum^n_{i=1} (\psi_i \otimes I_k(X)) (\psi_i \otimes
I_k(X))^*\right\|\cr
\le &\|\Psi_k(X,X^*)\|\cr
\le &\|\Psi\|_{cb}\|X\|^2,&(2.24)}$$
and a similar estimate holds for columns with $\psi_i$'s replaced by
$\theta_i$'s. It follows from (2.24) that
$$\|(\psi_1,\ldots, \psi_n)\|\le \|\Psi\|^{1/2}_{cb}\eqno (2.25)$$
and
$$\|(\theta_1,\ldots, \theta_n)^T\|\le \|\Theta\|^{1/2}_{cb}\eqno (2.26)$$
for any finite subset ${\cal F}$ of $\Lambda$. By [2, Theorem~3.1], $v\in
CB({\cal X}, {\cal N}) \otimes_{w^*h} CB({\cal Y}, {\cal N})$, and
$$\|v\|_{w^*h} \le \|\Psi\|^{1/2}_{cb} \|\Theta\|^{1/2}_{cb} \le 1+\vp.
\eqno (2.27)$$
Since $\vp>0$ was arbitrary, $\|v\|_{w^*h} \le |||v|||$ for all $v\in {\cal
V}$.

Conversely, suppose that $v\in CB({\cal X}, {\cal N}) \otimes_{w^*h}
CB({\cal Y}, {\cal N})$ and $\|v\|_{w^*h}=1$. By [2, Theorem~3.1], there
exists a representation $v = \sum\limits_{j\in\Lambda} \psi_j \otimes
\theta_j$, where the norms of the row of $\psi_j$'s and column of
$\theta_j$'s are both 1. For any finite subset
${\cal F} = \{1,2,\ldots, n\}$ of $\Lambda$ (after renumbering), let
$$\Psi_{\cal F}(x_1,x^*_2) = \sum^n_{j=1} \psi_j(x_1) \psi_j(x_2)^*\eqno
(2.28)$$
and
$$\Theta_{\cal F}(y^*_1,y_2) = \sum^n_{j=1}
\theta_j(y_1)^*\theta_j(y_2).\eqno (2.29)$$
It is then immediate from the definition of the norms in ${\Bbb
M}_k(CB({\cal X}, {\cal N}))$ and ${\Bbb M}_k(CB({\cal Y}, {\cal N}))$ that
$\|\Psi_{\cal F}\|_{cb}$, $\|\Theta_{\cal F}\|_{cb}\le 1$, from which it
follows that $\|\Psi\|_{cb}$, $\|\Theta\|_{cb}\le 1$. Thus $v\in {\cal V}$
and $|||v||| \le \|v\|_{w^*h}$, proving the reverse inequality.\medskip

 There is a natural map $\nu\colon \ CB(\ca X,\ca N)\otimes_{w^*h}
CB(\ca Y,\ca N)\to$  $CB^2(\ca X\times \ca Y,\ca N)$ defined by
$$\nu\left(\sum_{j\in\Lambda} \psi_j\otimes \theta_j\right)(x,y) =
\sum_{j\in\Lambda} \psi_j(x) \theta_j(y).\eqno (2.30)$$
The sum on the right hand side of (2.23) may be viewed as the product of
elements from $R(\Lambda) \overline\otimes \ca N$ and $C(\Lambda) \overline
\otimes \ca N$, from which the estimate
$$\eqalignno{\left\|\sum_{j\in\Lambda} \psi_j(x) \theta_j(y)\right\| &\le
\left\|\sum_{j\in\Lambda} \psi_j(x) \psi_j(x)^*\right\|^{1/2}
\left\|\sum_{j\in\Lambda} \theta_j(y)^* \theta_j(y)\right\|^{1/2}\cr
&= \|\Psi(x,x^*)\|^{1/2} \|\Theta(y^*,y)\|^{1/2}\cr
&\le \|\Psi\|^{1/2}_{cb} \|\Theta\|^{1/2}_{cb} \|x\|\, \|y\|&(2.31)}$$
is immediate. This inequality lifts to the $n$-fold amplification, showing
that $\nu$ is a contraction. Subsequently we show that $\nu$ is a complete
quotient map when $\ca N$ is injective.\vfill\eject

\n {\bf \S 3. ${\cal A}$-modular bilinear operators}

In this section we generalize Wittstock's one variable completely bounded
modular extension theorem [33] to two variables. The proof involves
standard techniques of modifying a completely bounded operator to one that
has a representation with good computational properties. Our approach would
apply to any number of variables, and in the case of one variable is
perhaps simpler than Wittstock's original method.

\proclaim Theorem 3.1. Let $\ca A$ be a $C^*$-subalgebra of $B(H)$, let
$\ca E$ and $\ca F$ be norm closed subspaces of $B(H)$ which are also $\ca
A$-modules, and let $\phi\colon \ \ca E\times \ca F\to B(H)$ be a
completely bounded $\ca A$-modular bilinear operator. Then $\phi$ extends
to a completely bounded $\ca A$-modular bilinear operator $\psi\colon \
B(H)\times B(H)\to B(H)$ with preservation of norm. Moreover, $\psi$ has a
representation
$$\psi(x,y) = V\pi(x) T\pi(y)W, \qquad x,y\in B(H),\eqno (3.1)$$
where
$\pi\colon \ B(H)\to B(K)$ is a representation, and $V,T,W$ are continuous
linear operators
$$H \brel W K \brel T K \brel V H\eqno (3.2)$$
satisfying $\|\psi\|_{cb} = \|V\|\, \|T\|\, \|W\|$, and
$$aV = V\pi(a),\quad \pi(a)T = T\pi(a),\quad \pi(a)W=W a.\eqno (3.3)$$
\medskip

\n {\bf Proof.} A completely bounded bilinear operator can be extended with
preservation of completely bounded norm, so let $\theta\colon \ B(H)\times
B(H) \to B(H)$ be any such extension of $\phi$. By (2.3), $\theta$ has a
representation
$$\theta(x,y) = V_1\pi(x) T_1\pi(y) W_1,\qquad x,y\in B(H)\eqno (3.4)$$
with $\|\theta\|_{cb} = \|V_1\|\, \|T_1\|\,\|W_1\|$. These operators are
successively replaced by inserting suitable projections from $B(K)$ into
(3.4).

By $\ca A$-modularity,
$$(aV_1-V_1\pi(a)) \pi(e) T_1\pi(f) W_1\xi = 0\eqno (3.5)$$
for all $a\in \ca A, e\in\ca E, f\in \ca F$, and $\xi\in H$. Let
$$K_1 = \overline{\rm span}\{\pi(e) T_1\pi(f) W_1\xi\colon \ e\in \ca E,
f\in \ca F, \xi\in H\}.$$
Then $K_1$ is a closed $\pi(A)$-invariant subspace of $K$, since $\pi(a)
\pi(e) = \pi(ae) \in \pi(\ca E)$. Thus the projection $P_1$ of $K$ onto
$K_1$ is in $\pi(\ca A)'$. Let $V=V_1P_1$ and define $\theta_1 \in
CB^2(B(H), B(H))$ by
$$\theta_1(x,y) = V\pi(x) T_1\pi(y)W_1,\qquad x,y\in B(H).\eqno (3.6)$$
Clearly $\|\theta_1\|_{cb} \le \|\theta\|_{cb}	= \|\phi\|_{cb}$ and we now
verify that $\theta_1$ is an extension of $\phi$. For all $\xi,\eta \in H$,
$e\in \ca E$, and $f\in \ca F$,
$$\eqalignno{\langle\theta_1(e,f)\xi,\eta\rangle &= \langle V_1P_1\pi(e)
T_1\pi(f) W_1\xi,\eta\rangle\cr
&= \langle V_1\pi(e) T_1\pi(f) W_1\xi, \eta\rangle\cr
&= \langle\theta(e,f)\xi,\eta\rangle\cr
&= \langle\phi(e,f)\xi,\eta\rangle,&(3.7)}$$
and so $\theta_1$ extends $\phi$. The second equality in (3.7) is immediate
from the definition of $K_1$. From (3.5), $aV_1-V_1\pi(a)$ annihilates
$K_1$ for all $a\in
\ca A$, so $(aV_1-V_1\pi(a))P_1=0$. Since $P_1$ commutes with $\pi(\ca A)$,
we obtain
$$aV-V\pi(a) = aV_1P_1- V_1P_1\pi(a) = (aV_1-V_1\pi(a))P_1=0 \eqno (3.8)$$
for all $a\in\ca A$. We now make a second modification to $\theta$.

By $\ca A$-modularity,
$$\eqalignno{ \langle V\pi(e) (\pi(a) T_1-T_1\pi(a))\pi(f)
W_1\xi,\eta\rangle
&= \langle(\theta_1(ea,f) - \theta_1(e,af))\xi,\eta\rangle\cr
&= \langle(\phi(ea,f) - \phi(e,af))\xi,\eta\rangle\cr
&= 0&(3.9)}$$
for all $a\in \ca A$, $e\in \ca E$, $f\in\ca F$ and $\xi,\eta \in H$. Thus
$$\langle (\pi(a)T_1-T_1\pi(a)) \pi(f)W_1\xi,\pi(e)^*
V^*\eta\rangle=0.\eqno (3.10)$$
Let
$$\eqalignno{K_2 &= \overline{\rm span}\{\pi(f) W_1\xi\colon \ f\in\ca F,
\xi\in H\},&(3.11)\cr
K_3 &= \overline{\rm span}\{\pi(e)^* V^*\eta\colon \ e\in \ca E, \eta \in
H\},&(3.12)}$$
and let $P_2$ and $P_3$ respectively be the projections onto these
subspaces of $K$. The relations
$$\pi(a) \pi(f) = \pi(af)\in \pi(\ca F), \quad \pi(a) \pi(e)^* =
\pi(ea^*)^* \in\pi(\ca E)^*\eqno (3.13)$$
show that $K_2$ and $K_3$ are invariant subspaces for $\pi(\ca A)$, and so
$P_2, P_3\in \pi(\ca A)'$. Let $T=P_3T_1P_2$, and define $\theta_2\in
CB^2(B(H), B(H))$ by
$$\theta_2(x,y) = V\pi(x) T\pi(y)W_1,\qquad x,y\in B(H).\eqno (3.14)$$
Then, for $e\in \ca E$, $f\in\ca F$, $\xi,\eta\in H$,
$$\eqalignno{\langle \theta_2(e,f)\xi,\eta\rangle &= \langle T_1P_2\pi(f)
W_1\xi, P_3\pi(e)^*V^*\eta\rangle\cr
&= \langle T_1\pi(f) W_1\xi, \pi(e)^*V^*\eta\rangle\cr
&= \langle \theta_1(e,f)\xi,\eta\rangle\cr
&= \langle \phi(e,f)\xi,\eta\rangle,&(3.15)}$$
verifying that $\theta_2$ extends $\phi$, and norm preservation is clear.
Moreover
$$\eqalignno{&\phantom{=} \langle(\pi(a) T_1-T_1\pi(a)) \pi(f) W_1\xi,
\pi(e)^* V^*\eta\rangle\cr
&= \langle(\theta_1(ea,f) - \theta_1(e,af))\xi,\eta\rangle\cr
&= \langle(\phi(ea,f) - \phi(e,af)) \xi,\eta\rangle\cr
&=0&(3.16)}$$
and so each operator $\pi(a) T_1-T_1\pi(a)$, $(a\in\ca A)$, maps
$K_2$ into $K^\bot_3$, by (3.11) and (3.12). Thus
$$\eqalignno{\pi(a) T-T\pi(a) &= \pi(a) P_3T_1P_2-P_3T_1P_2\pi(a)\cr
&= P_3(\pi(a)T_1-T_1\pi(a))P_2\cr
&= 0&(3.17)}$$
for $a\in\ca A$ since $P_2,P_2\in \pi(\ca A)'$. Thus (3.8) and (3.17) show
that, for $a\in \ca A$, $x,y\in B(H)$,
$$\theta_2(ax,y) = a\theta_2(x,y),\quad \theta_2(xa,y) =
\theta_2(x,ay).\eqno (3.18)$$

The final modification is to define $P_4$ to be the projection of $K$ onto
$$K_4  = \overline{\rm span}\{\pi(f)^*T^*\pi(e)^* V^*\eta\colon \ e\in\ca
E, f\in \ca F, \eta\in H\},\eqno (3.19)$$
and to let $\psi\in CB^2(B(H), B(H))$ be the bilinear operator obtained by
replacing $W_1$ by $W=P_4W_1$ in (3.14):
$$\psi(x,y) = V\pi(x) T\pi(y)W,\qquad x,y\in B(H).\eqno (3.20)$$
As before, one can verify that $\psi$ extends $\phi$ with preservation of
norm, and that the relations (3.2) hold. The details are so similar that we
omit them.
\medskip

\n {\bf Remark 3.2.}  Note that if $F=E^*$ and $\phi$ is in addition
completely positive, then $\psi$ can be chosen to be completely positive
with $W=V^*$ and $T\ge 0$.\medskip

We end this section with a technical result which will be needed
subsequently; since it deals with modular bilinear operators we include it
here. Below, the important point is the normality condition on the
representation.

\proclaim Proposition 3.3. Let $\ca R, \ca S\subseteq B(H)$ be a commuting
pair of von~Neumann algebras and let $\phi\in CB^2(C^*(\ca R,\ca S), B(H))$
be $\ca S$-modular. Then there exist a representation $\pi\colon \ C^*(\ca
R,\ca S)\to B(L)$, whose restriction to $\ca S$ is normal, and continuous
linear operators $W\colon \ H\to L$, $T\colon \ L\to L$ and $V\colon \ L\to
H$ such that
$$\phi(x,y) = V\pi(x) T\pi(y)W,\qquad x,y\in C^*(\ca R,\ca S),\eqno
(3.21)$$ $\|\phi\|_{cb} = \|V\|\, \|T\|\, \|W\|$, and
$$sV = V\pi(s),\quad \pi(s) T = T\pi(s),\quad \pi(s) W = Ws,\qquad
s\in \ca S.\eqno (3.22)$$

\n {\bf Proof.} By Theorem 3.1, $\phi$ has a representation
$$\phi(x,y) = V_1\rho(x) T_1\rho(y)W_1,\qquad x,y\in C^*(\ca R,\ca
S)\eqno (3.23)$$
where $\rho\colon \ C^*(\ca R,\ca S) \to B(K)$ is a representation and
$W_1\colon \ H\to K$, $T_1\colon \ K\to K$, and $V_1\colon \ K\to H$
satisfy $\|\phi\|_{cb} = \|V_1\|\, \|T_1\|\, \|W_1\|$. Moreover
$$sV_1 = V_1\rho(s),\quad \rho(s)T_1 = T_1\rho(s),\quad \rho(s)W_1=W_1s.
\eqno (3.24)$$
By the decomposition of a representation into its normal and singular parts
[30, Theorem~III.2.14] there is a central projection $p\in \rho(\ca S)''$
so that
$$s\to \langle \rho(s)p\xi,\eta\rangle,\quad s\to
\langle\rho(s)(1-p)\xi,\eta\rangle\eqno (3.25)$$
are respectively normal and singular linear functionals on $\ca S$ for all
$\xi,\eta\in K$. Now
$$\eqalignno{s\to \langle\phi(rs,y)\xi,\eta\rangle &=
\langle\phi(sr,y)\xi,\eta\rangle\cr
&= \langle s\phi(r,y)\xi,\eta\rangle\cr
&= \langle\rho(s) T_1\rho(y) W_1\xi, \rho(r)^*V^*_1\eta\rangle&(3.26)}$$
is normal on $\ca S$ for all $r\in \ca R$, $y\in C^*(\ca R,\ca S)$, and
$\xi,\eta\in H$ by the second equality in (3.26). Thus
$$\langle\phi(rs,y)\xi,\eta\rangle = \langle V_1\rho(rs) pT_1\rho(y)
W_1\xi,\eta\rangle\eqno (3.27)$$
for $r\in\ca R, s\in \ca S, y\in C^*(\ca R,\ca S)$, $\xi,\eta\in H$.

Now $p$ is central in $\rho(\ca S)''$ and so commutes with $\rho(\ca S)$.
Moreover, by hypothesis, $\rho(\ca R) \subseteq \rho(\ca S)' = (\rho(\ca
S)'')'$, so $p~(\in \rho(\ca S)'')$ commutes with $\rho(r)$ for $r\in
\ca R$. It follows that $p\in C^*(\ca R,\ca S)'$. By (3.24), $T_1\in
\rho(\ca S)'$. Thus $T_1$ commutes with all operators in $\rho(\ca S)''$,
and in particular with $p$. It follows from (3.27) that
$$\langle\phi(rs,y) \xi,\eta\rangle = \langle V_1p\rho(rs) p T_1p\rho(y)
pW_1\xi,\eta\rangle \eqno (3.28)$$
and hence, for $x,y\in C^*(\ca R,\ca S)$,
$$\phi(x,y) = V_1p\rho(x) pT_1p\rho(y)pW_1.\eqno (3.29)$$
Let $L\subseteq K$ be the range of $p$, and observe that this is invariant
for $C^*(\ca R,\ca S)$ since $p\in C^*(\ca R,\ca S)'$. Thus we may define
$\pi = \rho|_L$, $W = pW_1\colon \ H\to L$, $T = pT_1|_L\colon \ L\to L$,
and $V = V_1|_L\colon \ L\to H$, to obtain the representation
$$\phi(x,y) = V\pi(x) T\pi(y)W,\qquad x,y\in C^*(\ca R,\ca S)\eqno (3.30)$$
from (3.29). For $\xi,\eta\in L$ and $s\in\ca S$,
$$\langle\pi(s)\xi,\eta\rangle = \langle\pi(s) p\xi,\eta\rangle = \langle
\rho(s) p\xi,\eta\rangle,\eqno (3.31)$$
which, by (3.25), defines a normal functional on $\ca S$. Thus the
restriction of $\pi$ to $\ca S$ is a normal representation. The required
relations (3.22) follow easily from (3.24). For $\xi\in L$,
$$(sV-V\pi(s))\xi = (sV_1-V_1\rho(s))\xi = 0\eqno (3.32)$$
since $\rho(s) \xi\in L$, and so $sV = V\pi(s)$. The other two equalities
in (3.22) are verified similarly. Finally, $\|\phi\|_{cb} = \|V\|\, \|T\|\,
\|W\|$ is clear from the construction.\vfill\eject

\n {\bf \S 4. Factorization for injective ranges}

The first two propositions are steps on the way to Theorem 4.4, which is
the main result of this section. Throughout, the sum
$\sum\limits_{j\in\Lambda}$ of operators over the index set $\Lambda$ is
weakly convergent. Note that $\Lambda$ may be uncountable even for
von~Neumann algebras with separable predual if $\phi$ is not normal, for
example. The conclusions can be strengthened slightly, which is discussed
in Remark~4.5. In this section, all operators denoted by upper case greek
letters $\Psi, \Theta, \Gamma, \Delta$, are defined as in (2.20) and
(2.21).

We begin with a technical lemma which will be needed subsequently. The same
result for different pairs of tensor products may be found in [14, 17].
On the algebraic tensor product $\ca A\otimes \ca B \otimes \ca C\otimes
\ca D$ of $C^*$-algebras, the shuffle map $S$ is defined by
$$S(a\otimes b \otimes c \otimes d) = a\otimes c\otimes b \otimes d.\eqno
(4.1)$$

\proclaim Lemma 4.1. Let $\ca A, \ca B, \ca C$ and $\ca D$ be
$C^*$-algebras. The shuffle map induces a complete contraction
$$S\colon \ (\ca A\otimes_{\rm min} \ca B) \otimes_h (\ca C\otimes_{\rm min}
\ca D) \to (\ca A \otimes_h \ca C) \otimes_{\rm min} (\ca B\otimes_h \ca
D). \eqno (4.2)$$

\n {\bf Proof.} We will use the fact, [23], that there is an isometric
identification between\break $CB^2(\ca A\times \ca B, B(H))$ and $CB(\ca
A\otimes_h \ca B, B(H))$ for $C^*$-algebras $\ca A$ and $\ca B$. Let
\break $\phi_1\colon \ \ca A\otimes_h \ca C\to B(K_1)$ and $\phi_2\colon \
\ca B\otimes_h \ca D\to B(K_2)$ be completely isometric embeddings. Then
$\phi_1\otimes\phi_2\colon \ (\ca A\otimes_h \ca C) \otimes_{\rm min} (\ca
B\otimes_h \ca D) \to B(K_1) \otimes_{\rm min} B(K_2)$ is a completely
isometric embedding [1]. From (2.3) and (2.4), $\phi_1$ and $\phi_2$ may
be expressed by
$$\eqalignno{ \phi_1(a\otimes c) &= V_1\pi_1(a) T_1\rho_1(c)W_1,&(4.3)\cr
\phi_2(b\otimes d) &= V_2\pi_2(b) T_2\rho_2(d)W_2,&(4.4)}$$
where $\pi_i$, $\rho_i$ are $*$-representations and $V_i,T_i$ and $W_i$ are
contractive operators between appropriate Hilbert spaces, for $i=1,2$. By
the definition of the minimal tensor product [30], $\pi_1\otimes\pi_2$ and
$\rho_1\otimes\rho_2$ define $*$-representations of $\ca A\otimes_{\rm min}
\ca B$ and $\ca C\otimes_{\rm min} \ca D$ respectively. Define $\psi \in
CB^2((\ca A\otimes_{\rm min} \ca B) \times (\ca C\otimes_{\rm min} \ca
D)$, $B(K_1\otimes_2 K_2))$ by
$$\psi(a\otimes b, c\otimes d) = (V_1\otimes V_2) (\pi_1(a) \otimes
\pi_2(b)) (T_1\otimes T_2) (\rho_1(c) \otimes \rho_2(d))(W_1\otimes
W_2).\eqno (4.5)$$
Then $\psi$ is a completely contractive bilinear operator, and
$$\psi(a\otimes b, c\otimes d) = \phi_1(a\otimes c) \otimes \phi_2(b\otimes
d).\eqno (4.6)$$
Letting $\widetilde \psi$ be the associated completely contractive linear
operator on\hfil\break $(\ca A\otimes _{\rm min} \ca B) \otimes_h (\ca
C\otimes_{\rm
min}~\ca D)$, the complete contractivity of $S$ follows from the relation
$$S = (\phi_1\otimes \phi_2)^{-1}\widetilde \psi,\eqno (4.7)$$
where $(\phi_1\otimes\phi_2)^{-1}$ is defined on the range of
$\phi_1\otimes\phi_2$.

\proclaim Proposition 4.2. Let $\ca M\subseteq \ca N\subseteq B(H)$ be an
inclusion of injective von~Neumann algebras with $\ca M$ a factor. If
$\phi\in CB^2(\ca M,\ca N)$ then there exist $\psi_j,\theta_j \in CB(\ca
M,\ca N)$ satisfying
$$\phi(m_1,m_2) = \sum_{j\in\Lambda}\psi_j(m_1) \theta_j(m_2),\qquad
m_1,m_2\in \ca M\eqno (4.8)$$
and
$$\|\phi\|_{cb} = \|\Psi\|^{1/2}_{cb} \|\Theta\|^{1/2}_{cb}.\eqno (4.9)$$

\n {\bf Proof.} We will use the results of Effros and Lance [15,
Proposition~4.5 and Corollary~4.6] that for an injective von~Neumann algebra
$\ca M$ the map $\eta_{\ca M} \colon \ m\otimes m'\to mm'$ from $\ca M
\otimes_{\rm min}~\ca M'$ into $C^*(\ca M,\ca M')$ is a bounded surjective
$*$-homomorphism, and is additionally a $*$-isomorphism if $\ca M$ is a
factor. These were proved originally for semidiscrete von~Neumann algebras,
but semidiscreteness is equivalent to injectivity.

Define the operator $\tilde\phi\colon \ C^*(\ca M,\ca N') \times C^*(\ca M,
\ca N') \to C^*(\ca N,\ca N')$ by
$$\tilde\phi(m_1n'_1,m_2n'_2) = \phi(m_1,m_2) n'_1n_2'\eqno (4.10)$$
for $m_1,m_2\in \ca M$, $n'_1,n'_2\in \ca N'$. To see that $\tilde\phi$ is
a well defined completely bounded bilinear operator observe that it is the
composition of the following completely bounded maps.\medskip

\n (1)~~The inclusion
$$I\colon \ C^*(\ca M,\ca N')\to C^*(\ca M,\ca M').\eqno (4.11)$$
Here $\ca N'\subseteq \ca M'$ since $\ca M\subseteq \ca N$.

\n (2)~~The inverse of the Effros-Lance isomorphism
$$\eta^{-1}_{\ca M}\colon \ C^*(\ca M,\ca M')\to \ca M\otimes_{\rm min} \ca
M'\eqno (4.12)$$
which exists since $\ca M$ is an injective factor [15]. Note that the range
of $\eta^{-1}_{\ca M}I$ is $\ca M\otimes_{\rm min} \ca N'$.

\n (3)~~The completely bounded bilinear operator $\phi\otimes\lambda$ from
$(\ca M\otimes_{\rm min} \ca M') \times (\ca M\otimes_{\rm min} \ca M')$ to
$\ca N \otimes_{\rm min} \ca M'$ given by
$$(m_1\otimes m'_1) \times (m_2\otimes m'_2) \to \phi(m_1,m_2) \otimes
m'_1m'_2,\eqno (4.13)$$
where $\lambda\colon\ \ca M'\times \ca M'\to \ca M'$ is the completely
contractive bilinear multiplication map $\lambda(m'_1,m'_2) = m'_1m'_2$.
There are several ways to show that $\phi\otimes\lambda$ is completely
bounded. One such is to observe that the shuffle map $S\colon \ (\ca
A\otimes_{\rm min} \ca B) \otimes_h (\ca C\otimes_{\rm min}\ca D) \to (\ca
A\otimes_h \ca C) \otimes_{\rm min} (\ca B\otimes_h\ca D)$
 for
$C^*$-algebras $\ca A, \ca B, \ca C$ and $\ca D$, defined in (4.1), is
completely contractive, by Lemma 4.1.
 Then recall from [23] that a completely bounded
bilinear operator $\phi\colon \ \ca A\times \ca B\to \ca C$ induces a
completely bounded linear operator $\psi\colon \ \ca A\otimes_h \ca B\to
\ca C$ of the same norm by $a\otimes b\to \phi(a,b)$. Combining these
results, we obtain a completely bounded linear operator on $(\ca
M\otimes_{\rm min} \ca M') \otimes_h (\ca M\otimes_{\rm min} \ca M')$ by
$$(m_1\otimes m'_1) \otimes (m_2\otimes m'_2) \to (m_1 \otimes m_2) \otimes
(m'_1\otimes m'_2) \to \phi(m_1,m_2) \otimes m'_1m'_2\eqno (4.14)$$
and  this is the linearization of $\phi\otimes\lambda$. It is easy to check
that $\|\phi\otimes\lambda\|_{cb} = \|\phi\|_{cb}$. Note that $(\phi
\otimes \lambda)\eta^{-1}_{\ca M}I$ has range in $\ca M\otimes_{\rm min}
\ca N'$.

\n (4)~~The Effros-Lance homomorphism
$\eta_{\ca N}\colon\ \ca N\otimes_{\rm min} \ca N'\to C^*(\ca N,\ca N')$,
which is continuous since $\ca N$ is injective [15].

We now see that
$$\tilde\phi = \eta_{\ca N}(\phi\otimes \lambda) (\eta^{-1}_{\ca M} \otimes
\eta^{-1}_{\ca M})(I\otimes I)\eqno (4.15)$$
and so $\tilde\phi$ is well defined and completely bounded, as asserted, and
it is easily checked that $\|\tilde\phi\|_{cb} = \|\phi\|_{cb}$. Moreover
the construction of $\tilde\phi$ shows that this operator is $\ca
N'$-modular. By Proposition~3.3 (with $\ca R = \ca M, \ca S=\ca N'$),
$\tilde\phi$ has a representation
$$\tilde\phi(x,y) = V\pi(x) T\pi(y)W,\qquad x,y\in C^*(\ca M,\ca N')\eqno
(4.16)$$
where $\pi\colon \ C^*(\ca M,\ca N')\to B(L)$ is a representation whose
restriction to $\ca N'$ is normal,
$$n'V = V\pi(n'),\quad \pi(n')T = T\pi(n'),\quad \pi(n')W = Wn',\qquad
n'\in \ca N',\eqno (4.17)$$
and $\|\phi\|_{cb} = \|V\|\, \|T\|\, \|W\|$.
By the structure theory of normal representations [11] we may assume that
$$L \subseteq \ell_2(\Lambda,H) = \ell_2(\Lambda) \otimes_2 H,\eqno
(4.18)$$
where $\Lambda$ is a sufficiently large index set, $L$ is invariant for the
von~Neumann subalgebra $I\otimes \ca N'$ of $B(\ell_2(\Lambda))
\overline\otimes B(H)$, and $\pi(n') = (I\otimes n')q$ where $q$ is the
projection in $(I\otimes \ca N')' = B(\ell_2(\Lambda)) \overline\otimes \ca
N$ onto $L$. Writing matrices relative to the decomposition
$\ell_2(\Lambda) \otimes_2 H = L \oplus L^\bot$, we have
$$\tilde \phi(x,y) = (V,0) \left(\matrix{\pi(x)&0\cr 0&0\cr}\right)
\left(\matrix{T&0\cr 0&0\cr}\right) \left(\matrix{\pi(y)&0\cr
0&0\cr}\right) {W\choose 0}\eqno (4.19)$$
for all $x,y\in C^*(\ca M,\ca N')$. Since $T$ commutes with $\pi(\ca N')$
by (4.17), $\left(\matrix{T&0\cr 0&0\cr}\right)$ commutes with
$\left(\matrix{\pi(n')&0\cr 0&z\cr}\right)$ for all $z\in B(L^\bot)$, and
in particular $T$ commutes with $I\otimes n'$, taking $z = (1-q)(I\otimes n')
(1-q)$. Hence $\left(\matrix{T&0\cr 0&0\cr}\right) \in B(\ell_2(\Lambda))
\overline\otimes \ca N$. Similarly
$$\eqalignno{n'(V,0) &= (n'V,0)\cr
&= (V\pi(n'),0)\cr
&= (V,0) \left(\matrix{\pi(n')&0\cr 0&z\cr}\right)&(4.20)}$$
for all $z\in B(L^\bot)$ so the same choice for $z$ as before shows that
$n'(V,0) = (V,0)(I\otimes n')$ for all $n'\in \ca N'$. Thus $(V,0) \in
R(\Lambda)  \overline\otimes \ca N$ where $R(\Lambda)$ is the row space of
$B(\ell_2(\Lambda))$. A similar calculation  shows that ${W\choose 0} \in
C(\Lambda)\overline\otimes \ca N$ where $C(\Lambda)$ is the column space of
$B(\ell_2(\Lambda))$. Finally $\left(\matrix{\pi(m)&0\cr 0&0\cr}\right)$
commutes with $\left(\matrix{\pi(n')&0\cr 0&z\cr}\right)$ for all $m\in \ca
M$, $n'\in \ca N'$, $z\in B(L^\bot)$, so the same choice of $z$ shows that
$\left(\matrix{\pi(m)&0\cr 0&0\cr}\right)$ commutes with $I\otimes \ca N'$
and thus lies in $B(\ell_2(\Lambda)) \overline\otimes \ca N$.

Now for each $j\in \Lambda$, let $e_j$ be the orthogonal projection onto
the basis vector of $\ell_2(\Lambda)$ with a 1 in the $j^{\rm th}$ place
and 0 in every other position. Then define $\psi_j$, $\theta_j\colon \ \ca
M\to \ca N$ for each $j\in\Lambda$ by letting $\psi_j(m)$ be the $j^{\rm
th}$ component of $(V,0) \left(\matrix{\pi(m)&0\cr 0&0\cr}\right)$ in
$R(\Lambda) \overline\otimes \ca N$, and by letting $\theta_j(m)$ be the
$j^{\rm th}$ component of $\left(\matrix{T&0\cr 0&0\cr}\right)
\left(\matrix{\pi(m)&0\cr 0&0\cr}\right) \left(\matrix{W\cr
0\cr}\right)$ in $C(\Lambda) \overline\otimes \ca N$. Thus
$$\psi_j(m) = (V,0) \left(\matrix{\pi(m)&0\cr 0&0\cr}\right) (e_j\otimes
I)\eqno (4.21)$$
and
$$\theta_j(m) = (e_j\otimes I) \left(\matrix{T&0\cr 0&0\cr}\right)
\left(\matrix{\pi(m)&0\cr 0&0\cr}\right) {W\choose 0}\eqno (4.22)$$
for all $m\in \ca M$. From (4.21), (4.22) and the restriction of (4.19)
to $\ca M\times \ca M$, it follows that
$$\phi(m_1,m_2) = \sum_{j\in\Lambda} \psi_j(m_1) \theta_j(m_2)\eqno
(4.23)$$
for all $m_1,m_2\in\ca M$.

Now
$$\eqalignno{\Psi(m_1,m_2) &= \sum_{j\in\Lambda} \psi_j(m_1)
\psi_j(m^*_2)^*\cr
&= \sum_{j\in\Lambda} (V,0) \left(\matrix{\pi(m_1)&0\cr 0&0\cr}\right)
(e_j\otimes I) \left(\matrix{\pi(m_2)&0\cr 0&0\cr}\right) {V^*\choose 0}\cr
&= V\pi(m_1) I\pi(m_2)V^*&(4.24)}$$
which shows that $\Psi$ is completely bounded and that $\|\Psi\|_{cb}  =
\|V\|^2$. Similarly
$$\eqalignno{\Theta(m_1,m_2) &= \sum_{j\in\Lambda} \theta_j(m^*_1)^*
\theta_j(m_2)\cr
&= \sum_{j\in\Lambda} (W^*,0) \left(\matrix{\pi(m_1)&0\cr 0&0\cr}\right)
\left(\matrix{T^*&0\cr 0&0\cr}\right) (e_j\otimes I) \left(\matrix{T&0\cr
0&0\cr}\right) \left(\matrix{\pi(m_2)&0\cr 0&0\cr}\right) {W\choose 0}\cr
&= W^*\pi(m_1) T^*T\pi(m_2)W&(4.25)}$$
which shows that $\Theta$ is completely bounded and that $\|\Theta\|_{cb}
\le \|W\|^2 \|T\|^2$. Then (4.9) is an immediate consequence of
$\|\phi\|_{cb} \ge \|\Psi\|^{1/2}_{cb}\|\Theta\|^{1/2}_{cb}$ and
(2.22).\medskip

The next step is to remove the hypothesis $\ca M\subseteq \ca N$ from
Proposition~4.2.

\proclaim Proposition 4.3. Let $\ca M\subseteq B(H)$ be an injective factor
and let $\ca N$ be an injective von~Neumann algebra. If $\phi\in CB^2(\ca
M,\ca N)$ then there  exist $\psi_j,\theta_j \in CB(\ca M,\ca N)$
satisfying
$$\phi(m_1,m_2) = \sum_{j\in\Lambda} \psi_j(m_1) \theta_j(m_2), \qquad
m_1,m_2\in \ca M\eqno (4.26)$$
and
$$\|\phi\|_{cb} = \|\Psi\|^{1/2}_{cb} \|\Theta\|^{1/2}_{cb}.\eqno (4.27)$$

\n {\bf Proof.} Identify $\ca M$ with the subalgebra $I\otimes\ca M$ of
$\ca N \overline\otimes B(H)$, and let $e$ be a rank one projection in
$B(H)$. Then define $\tilde\phi\colon \ \ca M\times\ca M\to \ca
N\overline\otimes B(H)$ by $\tilde\phi(m_1,m_2) = \phi(m_1,m_2) \otimes e$,
and note that $\|\tilde\phi\|_{cb} = \|\phi\|_{cb}$. The injective factor
$\ca M$ is now a subalgebra of the injective von~Neumann algebra $\ca
N\overline\otimes B(H)$ so Proposition~4.2 applies to $\tilde\phi$. Hence
$$\tilde\phi(m_1,m_2) = \sum_{j\in\Lambda} \gamma_j(m_1) \delta_j(m_2),
\qquad m_1,m_2\in\ca M\eqno (4.28)$$
where $\gamma_j,\delta_j\in CB(\ca M,\ca N \overline\otimes B(H))$ and
$\|\tilde\phi\|_{cb} = \|\Gamma\|^{1/2}_{cb}  \|\Delta\|^{1/2}_{cb}$. Thus
$$\phi(m_1,m_2) \otimes e = \sum_{j\in\Lambda} (1\otimes e) \gamma_j(m_1)
\delta_j(m_2) (1\otimes e), \qquad m_1,m_2\in\ca M.\eqno (4.29)$$

Let $\{\xi_k\colon \ k\in\Omega\}$ be an orthonormal basis for $H$, let
$\{f_k\colon \ k\in\Omega\}$ be the associated rank one projections and
choose partial isometries $\{v_k\colon \ k\in\Omega\}$ such that $v^*_kv_k
= f_k$ and $v_kv^*_k=e$. Then define $\psi_{jk}$, $\theta_{jk}\in CB(\ca M,
\ca N)$ by
$$\psi_{jk}(m) \otimes v_k = (1\otimes e) \gamma_j(m) (1\otimes f_j)\eqno
(4.30)$$
and
$$\theta_{jk}(m) \otimes v^*_k = (1\otimes f_k) \delta_j(m) (1\otimes
e).\eqno (4.31)$$
It follows from (4.29) that
$$\phi(m_1,m_2) = \sum_{\scriptstyle j\in\Lambda\atop \scriptstyle
k\in\Omega} \psi_{jk}(m_1) \theta_{jk}(m_2),\eqno (4.32)$$
giving the required factorization. Moreover
$$\eqalignno{\Psi(m_1,m_2) \otimes e &= \sum_{\scriptstyle j\in\Lambda
\atop \scriptstyle k\in\Omega} (1\otimes e) \gamma_j(m_1) (1\otimes f_k)
\gamma_j(m^*_2)^* (1\otimes e)\cr
&= (1\otimes e) \Gamma(m_1,m_2)(1\otimes e)&(4.33)}$$
so $\|\Psi\|_{cb} \le \|\Gamma\|_{cb}$. Similarly $\|\Theta\|_{cb} \le
\|\Delta\|_{cb}$, and (4.27) follows from the corresponding result for
$\Gamma$ and $\Delta$.
\medskip

We come now to the main result of the section.

\proclaim Theorem 4.4. Let $\ca X$ and $\ca Y$ be operator spaces and let
$\ca N\subseteq B(K)$ be an injective von~Neumann algebra. If $\phi\in
CB^2(\ca X\times \ca Y,\ca N)$ then there exist $\psi_j\in CB(\ca X,\ca N)$
and $\theta_j \in CB(\ca Y,\ca N)$ such that
$$\phi(x,y) = \sum_{j\in\Lambda} \psi_j(x) \theta_j(y),\qquad x\in \ca
X, y\in \ca Y,\eqno (4.34)$$
and
$$\|\phi\|_{cb} = \|\Psi\|^{1/2}_{cb} \|\Theta\|^{1/2}_{cb}.\eqno (4.35)$$

\n {\bf Proof.} Let $H$ be a Hilbert space such that $\ca X, \ca Y\subseteq
B(H)$
as operator spaces, and extend $\phi$ to $\phi_1\colon \ B(H) \times B(H)
\to B(K)$ with preservation of norm [23]. If $\bb E$ is the conditional
expectation from $B(K)$ to $\ca N$, let $\phi_2 = \bb E\phi_1$. Then
$\phi_2$ is a norm preserving extension of $\phi$ and $\phi_2\in CB^2(B(H),
\ca N)$. Proposition~4.3 now gives the required completely bounded linear
operators by restricting to $\ca X$ and  $\ca Y$ those defined on $B(H)$.
\medskip

\n {\bf Remark 4.5.} (i)~~We have stated Theorem~4.4 in the full generality
of operator spaces so that generalizations to multilinear operators are
immediate. For example, if \break $\phi\colon \ \ca A\times \ca B \times
\ca C \to
\ca N$ were a completely bounded trilinear operator, a factorization could
be obtained by identifying $\phi$ with a bilinear operator $\psi\in
CB^2((\ca A\otimes_h\ca B)\times \ca C,\ca N)$ and applying Theorem~4.4.
Even if $\ca A$ and $\ca B$ were $C^*$-algebras, $\ca A\otimes_h \ca B$ is
only an operator space is general. This technique of employing the Haagerup
tensor product comes from [23].

\n (ii)~~An examination of the proofs of the first two propositions shows
that if $\phi$ is separately normal in each variable, then the resulting
operators $\psi_j,\theta_j$ may be chosen to be normal. However this is not
necessarily true in Theorem~4.4 for dual operator spaces since the
conditional expectation may destroy normality.

\n (iii)~~If $\ca Y = \ca X^*$ in Theorem~4.4 then it makes sense to
consider the extra hypothesis of complete positivity for $\phi$. Again an
examination of the proofs, being careful to choose completely positive
extensions at each stage, shows that we may take $\theta_j$ to be
$\psi^*_j$, where $\psi^*_j(x^*)$ is defined to be $\psi_j(x)^*$.

\n (iv)~~For the simplest injective von Neumann algebra ${\Bbb C}$,
Theorem~4.4 recaptures, in different language, the result of [2] that the
dual of $\ca X\otimes_h \ca Y$ is $\ca X^* \otimes_{w^*h}\ca Y^*$ for
operator spaces $\ca X$ and $\ca Y$. Thus Theorem~4.4 may be viewed as a
generalization of [2, Theorem~3.2].

\vfill\eject

\n {\bf\S 5.  Completely bounded factorization implies injectivity}

In the previous section we considered factorizations $\phi(x,y) =
\sum\limits_{j\in\Lambda} \psi_j(x) \theta_j(y)$ in\break $CB^2(\ca X\times
\ca Y,\ca
N)$ of completely bounded bilinear operators, where the associated\break
 bilinear
operators $\Psi$ and $\Theta$ were completely bounded. We now wish to
broaden the set of admissible factorizations by considering ones for which
there exists a constant $K$ such that
$$\eqalignno{\left\|\sum_{j\in\Lambda} \psi_j(x) \psi_j(x)^*\right\| &\le
K\|x\|^2,&(5.1)\cr
\left\|\sum_{j\in\Lambda} \theta_j(y)^* \theta_j(y)\right\| &\le
K\|y\|^2.&(5.2)}$$
Then the associated bilinear operators $\Psi$ and $\Theta$ are still
completely positive and bounded, but perhaps not completely bounded. We
will distinguish these two factorizations by calling them type CB (for
completely bounded) and type B (for bounded) respectively. We emphasize
that the operators $\psi_j$, $\theta_j$ are completely bounded in both
cases, and the names reflect the nature of $\Psi$ and $\Theta$.

The following lemma records a standard decomposition of certain completely
bounded bilinear operators as a linear combination of continuous completely
positive bilinear operators in exactly the correct form for subsequent use.
The technique is well known in the theory of quadratic forms.

\proclaim Lemma 5.1. Let $\ca M$ and $\ca N$ be von~Neumann algebras and
suppose that $\phi \in CB^2(\ca M,\ca N)$ has a factorization of type B.
Then there exist continuous completely positive bilinear maps $\phi_k\colon
\ \ca M\times \ca M\to \ca N$, $1\le k \le 4$, such that
$$\phi = (\phi_1-\phi_2) + i(\phi_3-\phi_4).\eqno (5.3)$$

\n {\bf Proof.} Suppose that
$$\phi(m_1,m_2) = \sum_{j\in\Lambda} \psi_j(m_1) \theta_j(m_2),\qquad
m_1,m_2\in\ca M\eqno (5.4)$$
is a factorization with $\psi_j,\theta_j\in CB(\ca M,\ca N)$ satisfying
(5.1) and (5.2). The algebraic identity
$$\psi_j(m_1) \theta_j(m_2) = {1\over 4} \sum^4_{k=1} i^k(\psi_j(m_1) +
i^k\theta^*_j(m_1)) (\psi^*_j(m_2) + i^{-k}\theta_j(m_2)) \eqno (5.5)$$
expresses each bilinear operator $\psi_j(m_1) \theta_j(m_2)$ as a linear
combination of continuous completely positive bilinear operators. then
(5.3) follows by summing (5.5) over $j\in\Lambda$, provided that the
resulting sums on the right hand side define continuous bilinear operators.
We examine $\sum\limits_{j\in\Lambda} (\psi_j(m_1) + \theta^*_j(m_1))
(\psi^*_j(m_2) +\theta_j(m_2))$, which is typical. If $\xi,\eta$ are
arbitrary vectors, then
$$\eqalignno{&\phantom{\le} \left|\sum_{j\in\Lambda} \langle(\psi_j(m_1) +
\theta^*_j(m_1)) (\psi^*_j(m_2) + \theta_j(m_2)) \xi,\eta\rangle\right|\cr
&\le \sum_{j\in\Lambda} |\langle\psi_j(m^*_2)^*
\xi,\psi_j(m_1)^*\eta\rangle | + \sum_{j\in\Lambda} |\langle \theta_j(m_2)
\xi,\theta_j(m^*_1)\eta\rangle|\cr
&\quad + \sum_{j\in\Lambda} |\langle\theta_j(m_2)\xi,
\psi_j(m_1)^*\eta\rangle| + \sum_{j\in\Lambda} |\langle\psi_j(m^*_2)^* \xi,
\theta_j(m^*_1)\eta\rangle|.&(5.6)}$$
The estimation of each of these four terms is identical; we take the first
as typical. Then, applying the Cauchy-Schwarz inequality,
$$\eqalignno{&\phantom{\le} \sum_{j\in \Lambda} |\langle\psi_j(m^*_2)^*
\xi, \psi_j(m_1)^*\eta\rangle|\cr
&\le \sum_{j\in\Lambda} \|\psi_j(m^*_2)^* \xi\|\, \|\psi_j(m_1)^*\eta\|\cr
&\le \left(\sum_{j\in\Lambda} \|\psi_j(m^*_2)^* \xi\|^2\right)^{1/2}
\left(\sum_{j\in\Lambda} \|\psi_j(m_1)^* \eta\|^2\right)^{1/2}\cr
&= \left(\sum_{j\in\Lambda} \langle\psi_j(m^*_2)\psi_j(m^*_2)^*\xi,\xi
\rangle\right)^{1/2} \left(\sum_{j\in\Lambda} \langle\psi_j(m_1)
\psi_j(m_1)^* \eta,\eta\rangle\right)^{1/2}\cr
&\le \left\|\sum_{j\in\Lambda} \psi_j(m^*_2) \psi_j(m^*_2)^*\right\|^{1/2}
\left\|\sum_{j\in\Lambda} \psi_j(m_1) \psi_j(m_1)^*\right\|^{1/2} \|\xi\|\,
\|\eta\|\cr
&\le K\|m_1\|\, \|m_2\|\, \|\xi\|\, \|\eta\|&(5.7)}$$
by (5.1) and (5.2). Returning to (5.6), we obtain
$$\left\|{1\over 4} \sum_{j\in\Lambda} (\psi_j(m_1) + \theta^*_j(m_1))
(\psi^*_j(m_2) + \theta_j(m_2))\right\| \le K\|m_1\|\, \|m_2\|, \eqno
(5.8)$$
and so each sum over $j\in\Lambda$ on the right hand side of (5.5) is a
continuous bilinear operator of norm at most $K$.\medskip

Every infinite dimensional von Neumann algebra $\ca M$ contains a copy of
$\ell^\infty$. Let us fix such a copy, and denote by $\ca U$ its abelian
(and hence amenable) unitary group. Then let $\beta$ be a fixed normalized
invariant mean on the space of complex valued bounded functions on $\ca U$.
Letting $B^2(\ca M,\ca N)$ denote the space of bounded bilinear maps on
$\ca M\times \ca M$ into a von~Neumann algebra $\ca N$, there is an induced
map $\gamma\colon
\ B^2(\ca M,\ca N)\to B^2(\ca M,\ca N)$ defined as follows. For $x,y\in \ca
M$ and $\omega\in \ca N_*$ the function $f_{x,y,\omega}(u) =
\omega(\phi(xu,u^*y))$ is
bounded by $\|\omega\|\, \|\phi\|\, \|x\|\, \|y\|$. We define
$\gamma\phi(x,y) \in (\ca N_*)^* = \ca N$ by
$$\gamma\phi(x,y)(\omega) = \beta(f_{x,y,\omega}).\eqno (5.9)$$
The technique of averaging between the variables in the next lemma comes
from [12].

\proclaim Lemma 5.2. The map $\gamma\colon \ B^2(\ca M,\ca N)\to B^2(\ca
M,\ca N)$ is a linear contraction and $\gamma\phi$ satisfies
$$\gamma\phi(xa,y) = \gamma\phi(x,ay)\qquad x,y\in\ca M,\quad
a\in\ell^\infty \eqno (5.10)$$
for $\phi\in B^2(\ca M,\ca N)$. Moreover $\gamma$ preserves both complete
boundedness and complete positivity.

\n {\bf Proof.} Equation (5.10) is a consequence of the invariance of
$\beta$ and the standard fact that every unital $C^*$-algebra is the span
of its unitary group. The other parts of the lemma are routine deductions
from the definition of $\gamma$. For example, if $\ca N$ is represented on
$H$, $\xi_1,\ldots, \xi_n\in H$, $(x_{ij}) \in \bb M_n(\ca M)$, and $\phi$
is completely positive, then
$$\left\langle (\gamma\phi)_n ((x_{ij}), (x_{ij})^*) \left(\matrix{\xi_1\cr
\vdots\cr \xi_n\cr}\right), \left(\matrix{\xi_1\cr \vdots\cr
\xi_n\cr}\right) \right\rangle \ge 0\eqno (5.11)$$
because this inner product is obtained by applying $\beta$ to the
non-negative function
$$u\to \left\langle \phi_n((x_{ij}u), (x_{ij}u)^*)
\left(\matrix{\xi_1\cr \vdots\cr \xi_n\cr}\right), \left(\matrix{\xi_1\cr
\vdots\cr \xi_n\cr}\right)\right\rangle.$$
This shows that $\gamma\phi$ is
completely positive, and we omit further details.

\proclaim Theorem 5.3. Let $\ca M$ be an infinite dimensional von~Neumann
algebra and let $\ca N$ be a von~Neumann algebra. If every $\phi\in
CB^2(\ca M,\ca N)$ is a linear combination of continuous completely
positive bilinear operators, then $\ca N$ is injective.

\n {\bf Proof.} This result is deduced from a theorem of Haagerup [20,
Theorem~2.1] on decomposable completely bounded linear operators from
$\ell^\infty$ into a von~Neumann algebra $\ca N$. Fix a copy of
$\ell^\infty$ with unitary group $\ca U$ in $\ca M$ and an invariant mean
$\beta$ as in the previous lemma. Let $\bb E$ be the conditional
expectation from $\ca M$ onto $\ell^\infty$. If $\phi\in CB(\ell^\infty,\ca
N)$, define $\theta\in CB^2(\ca M,\ca N)$ by
$$\theta(m_1,m_2) = \phi(\bb E(m_1m_2)),\qquad m_1,m_2\in\ca M.\eqno
(5.12)$$
Since $m_1\times m_2\to m_1m_2$ is a completely bounded bilinear operator
on $\ca M\times \ca M$ and $\bb E$ is completely positive, it is routine to
check that $\theta$ is completely bounded. By hypothesis, there are
continuous completely positive bilinear operators $\theta_j\colon \ \ca
M\times \ca M\to \ca N$ $(j=1,\ldots, 4)$ such that
$$\theta = \theta_1-\theta_2 + i\theta_3-i\theta_4.\eqno (5.13)$$
Now $\theta(m_1u,u^*m_2) = \theta(m_1,m_2)$ for all $u\in\ca U$, by (5.12),
and so $\gamma\theta =\theta$. Applying $\gamma$ to (5.13) and noting that
$\phi(x) = \theta(x,1)$ for $x\in\ell^\infty$, we obtain
$$\phi(x) = \gamma\theta_1(x,1) - \gamma\theta_2(x,1) +
i\gamma\theta_3(x,1) - i\gamma\theta_4(x,1),\qquad x\in\ell^\infty. \eqno
(5.14)$$
Now define $\phi_j\colon \ \ell^\infty\to \ca N$ $(j=1,\ldots, 4)$ by
$\phi_j(x) = \gamma\theta_j(x,1)$. For $x\in \ell^\infty$,
$$\phi_j(xx^*) = \gamma\theta_j(xx^*,1) = \gamma\theta_j(x,x^*)\ge 0,\eqno
(5.15)$$
by Lemma~5.2 and complete positivity of $\theta_j$. Positivity of $\phi_j$
is immediate from (5.15), and since $\ell^\infty$
is an abelian $C^*$-algebra, complete positivity of $\phi_j$ follows from
[29]. Thus (5.14) expresses $\phi$ as a linear combination of completely
positive maps, and the result follows from [20, Theorem~2.1].
\medskip

Combining Lemma 5.1 and Theorem~5.3, we immediately have the following
result.

\proclaim Theorem 5.4. Let $\ca M$ and $\ca N$ be von~Neumann algebras with
$\ca M$ infinite dimensional. If every $\phi \in CB^2(\ca M,\ca N)$ has a
factorization of type $B$, then $\ca N$ is injective.\vfill\eject

\n {\bf \S 6. Summary of results}

In this section we collect together the results of previous sections.
Recall that a factorization $\phi(x,y) = \sum\limits_{j\in\Lambda}
\psi_j(x) \theta_j(y)$ is of type CB (respectively type B) if the
associated bilinear maps $\Psi$ and $\Theta$ are completely bounded
(respectively bounded). Also recall the product map $\nu\colon \ CB(\ca
X,\ca M) \otimes_{w^*h} CB(\ca Y,\ca M)\to CB^2(\ca X\times \ca Y,\ca M)$
from Section~2.

\proclaim Theorem 6.1. Let $\ca M$ be a von~Neumann algebra. Then the
following are equivalent:
\item{(i)} $\ca M$ is injective,
\item{(ii)} each $\phi\in CB^2(\ca M,\ca M)$ has a type CB factorization
$$\phi(m_1,m_2) = \sum_{j\in\Lambda} \psi_j(m_1) \theta_j(m_2),\qquad
\psi_j, \theta_j\in CB(\ca M,\ca M),$$
\item{(iii)} each $\phi\in CB^2(\ca M,\ca M)$ has a type $B$ factorization
with $\psi_j$, $\theta_j\in CB(\ca M,\ca M)$,
\item{(iv)} $\nu\colon \ CB(\ca M,\ca M) \otimes_{w^*h} CB(\ca M,\ca M)\to
CB^2(\ca M,\ca M)$ is surjective and a complete quotient map.\medskip

\n {\bf Proof.} (i) $\Rightarrow$ (ii). This is Theorem 4.4 in the case
$\ca X = \ca Y = \ca N = \ca M$.

\n (ii) $\Rightarrow$ (iii). This is obvious.

\n (iii) $\Rightarrow$ (i). This is Theorem~5.4. Of course there is nothing
to prove if $\ca M$ is finite dimensional.

\n (iv) $\Rightarrow$ (ii). This is immediate from the definition of the
$w^*$-Haagerup tensor product (see Section~2).

\n (i) $\Rightarrow$ (iv). The surjectivity of $\nu$ is Theorem~4.4. The
fact that $\nu$ is a complete quotient map is just a restatement of (4.35)
of that theorem.\medskip

\n {\bf Remark 6.2.} In Theorem 6.1 we have only given the main
equivalences that are internal to $\ca M$, but there are many others which
could be extracted from the previous sections. For example (ii) and (iii)
could be recast for general operator spaces $\ca X$ and $\ca Y$, while
dropping the complete quotient map hypothesis from (iv) gives a statement
which is clearly equivalent to (ii). \vfill\eject

\centerline{\bf References}

\item{ 1. } D.P.\ Blecher and V.I.\ Paulsen, Tensor products of operator
spaces, J.\ Funct.\ Anal., 99 (1991), 262--292.

\item{ 2. } D.P.\ Blecher and R.R.\ Smith, The dual of the Haagerup tensor
product, J.\ London Math.\ Soc., 45 (1992), 126--144.

\item{ 3. } J.W.\ Bunce and W.L.\ Paschke, Quasi-expectations and amenable
von~Neumann algebras, Proc.\ Amer.\ Math.\ Soc., 71 (1978), 232--236.

\item{ 4. } E.\ Christensen, F.\ Pop, A.M.\ Sinclair and R.R.\ Smith, On
the cohomology groups of certain finite von~Neumann algebras, Math.\ Ann.,
307 (1997), 71--92.

\item{ 5. } E.\ Christensen and A.M.\ Sinclair, Representations of
completely bounded multilinear operators, J.\ Funct.\ Anal., 72 (1987),
151--181.

\item{ 6. } E.\ Christensen and A.M.\ Sinclair, A survey of completely
bounded operators, Bull.\ London Math.\ Soc., 21 (1989), 417--448.

\item{ 7. } E.\ Christensen and A.M.\ Sinclair, On von~Neumann algebras
which are complemented subspaces of $B(H)$, J.\ Funct.\ Anal., 122 (1994),
91--102.

\item{ 8. } E.\ Christensen and A.M.\ Sinclair,  Module mappings into
von~Neumann algebras and injectivity, Proc.\ London Math.\ Soc., 71 (1995),
618--640.

\item{ 9. } E.\ Christensen and A.M.\ Sinclair, A cohomological
characterization of approximately finite dimensional von~Neumann algebras,
preprint.

\item{ 10. } A.\ Connes, Classification of injective factors, Ann.\ Math.,
104 (1976), 73--115.

\item{ 11. } J.\ Dixmier, Les alg\`ebres d'op\'erateurs dans l'espace
Hilbertien, Gauthier-Villars, Paris, 1969.

\item{ 12. } E.G.\ Effros, Amenability and virtual diagonals for
von~Neumann algebras, J.\ Funct.\ Anal., 78 (1988), 137--153.

\item{ 13. } E.G.\ Effros and A.\ Kishimoto, Module maps and
Hochschild-Johnson cohomology, Indiana Math.\ J., 36 (1987), 257--276.

\item{ 14. } E.G.\ Effros, J.\ Kraus and Z.J.\ Ruan, On two quantized
tensor products, in ``Operator algebras, mathematical physics and
low-dimensional topology (Istanbul 1991)'', pp.~125--145, Research Notes
Math. 5, Peters, Wellesley, Mass., 1993.

\item{ 15. } E.G.\ Effros and E.C.\ Lance, Tensor products of operator
algebras, Adv.\ in Math., 25 (1977), 1--34.

\item{ 16. } E.G.\ Effros and Z.J.\ Ruan, On approximation properties for
operator spaces, Internat.\ J.\ Math., 1 (1990), 163--187.

\item{ 17. } E.G.\ Effros and Z.J.\ Ruan, Operator convolution algebras:\
an approach to quantum groups, preprint.

\item{ 18. } U.\ Haagerup, The $\alpha$-tensor product for $C^*$-algebras,
unpublished manuscript, 1980.

\item{ 19. } U.\ Haagerup, A new proof of the equivalence of injectivity
and hyperfiniteness for factors  on a separable Hilbert space, J.\ Funct.\
Anal., 62 (1987), 160--201.

\item{ 20. } U.\ Haagerup, Injectivity and decomposition of completely
bounded
maps, Springer Lecture Notes in Math., 1132, Springer-Verlag, 1985,
pp.~170--222.

\item{ 21. } U.\ Haagerup and G.\ Pisier, Bounded linear operators between
$C^*$-algebras, Duke Math.\ J., 71 (1993), 889--925.

\item{ 22. } V.I.\ Paulsen, Completely bounded maps and dilations, Notes in
Mathematics Series 146, Pitman, New York, 1986.

\item{ 23. } V.I.\ Paulsen and R.R.\ Smith, Multilinear maps and tensor
norms on operator systems, J.\ Funct.\ Anal., 73 (1987), 258--276.

\item{ 24. } G.\ Pisier, Remarks on complemented subspaces of von~Neumann
algebras, Proc.\ Roy.\ Soc.\ Edinburgh, 121A (1991), 1--4.

\item{ 25. } G.\ Pisier, Projections from a von~Neumann algebra onto a
subalgebra, Bull.\ Soc.\ Math.\ France, 123 (1995), 139--153.

\item{ 26. } S.\ Popa, A short proof of ``injectivity implies
hyperfiniteness'' for finite von~Neumann algebras, J.\ Operator Theory, 16
(1986), 261--272.

\item{ 27. } A.M.\ Sinclair and R.R.\ Smith, Hochschild cohomology of
von~Neumann algebras, London Math.\ Soc.\ Lecture Note Series 203,
Cambridge Univ.\ Press, Cambridge, 1995.

\item{ 28. } R.R.\ Smith, Completely bounded module maps and the Haagerup
tensor product, J.\ Funct.\ Anal., 102 (1991), 156--175.

\item{ 29. } W.F.\ Stinespring, Positive functions on $C^*$-algebras,
Proc.\ Amer.\ Math.\ Soc., 6 (1955), 211--216.

\item{ 30. } M.\ Takesaki, Theory of operator algebras I, Springer-Verlag,
Berlin, 1979.

\item{ 31. } J.\ Tomiyama, On the projections of norm one in $W^*$-algebras
I, Proc.\ Japan Acad., 33 (1957), 608--612.

\item{ 32. } S.\ Wassermann, Injective $W^*$-algebras, Math.\ Proc.\ Cam.\
Phil.\ Soc., 82 (1977), 39--47.

\item{ 33. } G.\ Wittstock, Extensions of completely bounded
$C^*$-module homomorphisms, in Proc.\ Conference on Operator Algebras and
Group Representations, Neptun 1980, Pitman, New York, 1983.

\bye